\documentclass[11pt,reqno]{amsart}
\usepackage{amsmath}
\usepackage{amssymb}
\usepackage{amstext}
\usepackage{amsfonts}
\usepackage{amsthm}
\usepackage{ifthen}
\usepackage{xspace}
\usepackage{comment}
\usepackage{graphicx}
\usepackage{mathrsfs}
\usepackage{enumerate}
\usepackage{breakcites}

\RequirePackage[colorlinks,citecolor=blue,urlcolor=blue]{hyperref}
\RequirePackage{hypernat}

\numberwithin{equation}{section}  

\theoremstyle{plain}
\newtheorem{theorem}{Theorem}[section]
\newtheorem{corollary}[theorem]{Corollary}
\newtheorem{lemma}[theorem]{Lemma}

\theoremstyle{definition}

\theoremstyle{plain}
\newtheorem*{theorem*}{Theorem}
\newtheorem*{corollary*}{Corollary}
\newtheorem*{lemma*}{Lemma}
\newtheorem*{proposition*}{Proposition}

\theoremstyle{definition}
\newtheorem*{remark*}{Remark}
\newtheorem*{remarks*}{Remarks}
\newtheorem*{example*}{Example}
\newtheorem*{examples*}{Examples}
\newtheorem*{definition*}{Definition}
\newtheorem*{conjecture*}{Conjecture}
\newtheorem*{assumption*}{Assumption}
\newtheorem*{assumptions*}{Assumptions}

\def\RR{\mathbb{R}}

\def\EE{\mathbb{E}}

\begin{document}

\title{On the problem of reversibility of the entropy power inequality}

\dedicatory{Dedicated to Friedrich G\"otze on the occasion of his sixtieth birthday}

\author{Sergey G. Bobkov}
\address{Sergey G. Bobkov, School of Mathematics, University of Minnesota,
Vincent Hall 228, 206 Church St SE, Minneapolis MN 55455, USA}
\email{bobkov@math.umn.edu}
\thanks{Sergey G. Bobkov was supported in part by the NSF grant DMS-1106530.}

\author{Mokshay M. Madiman}
\address{Mokshay M. Madiman, Department of Statistics, Yale University,
24 Hillhouse Avenue, New Haven CT 06511, USA}
\email{mokshay.madiman@yale.edu}
\thanks{Mokshay M. Madiman was supported in part by the NSF CAREER grant DMS-1056996.}


\begin{abstract}
As was shown recently by the authors, the entropy power inequality can be 
reversed for independent summands with sufficiently concave densities,
when the distributions of the summands are put in a special position.
In this note it is proved that reversibility is impossible over the whole
class of convex probability distributions. Related phenomena for
identically distributed summands are also discussed.
\end{abstract}

\maketitle

\section{The reversibility problem for the entropy power inequality}
\label{bm-sec:intro}

Given a random vector $X$ in $\RR^n$ with density $f$, introduce the 
entropy functional (or Shannon's entropy)
$$
h(X) = - \int_{\RR^n} f(x) \log f(x)\,dx,
$$
and the entropy power
$$
H(X) = e^{2h(X)/n},
$$
provided that the integral exists in the Lebesgue sense. For example, 
if $X$ is uniformly distributed in a convex body $A \subset \RR^n$, 
we have
$$
h(X) = \log |A|, \qquad H(X) =|A|^{2/n},
$$
where $|A|$ stands for the $n$-dimensional volume of $A$.

The entropy power inequality due to Shannon and Stam indicates that
\begin{equation}\label{bm-epi}
H(X+Y) \geq H(X) + H(Y),
\end{equation}
for any two independent random vectors $X$ and $Y$ in $\RR^n$, for which
the entropy is defined (\cite{bm-Sha, bm-Sta}, cf. also \cite{bm-CC, bm-DCT, bm-SV}). 
This is one of the fundamental results in Information 
Theory, and it is of large interest to see how sharp \eqref{bm-epi} is.

The equality here is only achieved, when $X$ and $Y$ have normal
distributions with proportional covariance matrices. Note that the 
right-hand side is unchanged when $X$ and $Y$ are replaced with affine
volume-preserving transformation, that is, with random vectors
\begin{equation}\label{bm-volpres}
\widetilde X = T_1(X), \qquad \widetilde Y = T_2(Y) \qquad
(|{\rm det} T_1| = |{\rm det} T_2| = 1).
\end{equation}
On the other hand, the entropy power $H(\widetilde X + \widetilde Y)$
essentially depends on the choice of $T_1$ and $T_2$. Hence, it is 
reasonable to consider a formally improved variant of \eqref{bm-epi},
\begin{equation}\label{bm-epi2}
\inf_{T_1,T_2} H(\widetilde X + \widetilde Y) \geq H(X) + H(Y),
\end{equation}
where the infimum is running over all affine maps 
$T_1,T_2:\RR^n \rightarrow \RR^n$ subject to \eqref{bm-volpres}. 
(Note that one of these maps may be taken to be the identity operator.) 
Now, equality in \eqref{bm-epi2} is achieved, whenever $X$ and $Y$ have normal
distributions with arbitrary positive definite covariance matrices.

A natural question arises: When are both the sides of \eqref{bm-epi2} of a similar 
order? For example, within a given class of probability distributions 
(of $X$ and $Y$), one wonders whether or not it is possible to reverse 
\eqref{bm-epi2} to get
\begin{equation}\label{bm-repi-q}
\inf_{T_1,T_2} H(\widetilde X + \widetilde Y) \leq C(H(X) + H(Y))
\end{equation}
with some constant $C$. 

The question is highly non-trivial already for the class of uniform 
distributions on convex bodies, when it becomes to be equivalent 
(with a different constant) to the inverse Brunn-Minkowski inequality
\begin{equation}\label{bm-rBM}
\inf_{T_1,T_2} \big|\widetilde A + \widetilde B\big|^{1/n} \leq C 
\left(|A|^{1/n} + |B|^{1/n}\right).
\end{equation}
Here 
$\widetilde A + \widetilde B = \{x+y: x \in \widetilde A, \ y \in \widetilde B\}$
stands for the Minkowski sum of the images
$\widetilde A = T_1(A)$, $\widetilde B = T_2(B)$ of arbitrary convex bodies 
$A$ and $B$ in $\RR^n$. To recover such an equivalence, one takes for 
$X$ and $Y$ independent random vectors uniformly distributed in $A$ and $B$. 
Although the distribution of $X+Y$ is not uniform in $A+B$,
there is a general entropy-volume relation
$$
\frac{1}{4}\
|A + B|^{2/n} \leq\, H(X+Y) \leq\ |A + B|^{2/n},
$$
which may also be applied to the images $\widetilde A, \widetilde B$
and $\widetilde X$, $\widetilde Y$ (cf. \cite{bm-BM3}).

The inverse Brunn-Minkowski inequality \eqref{bm-rBM} is indeed true and represents 
a deep result in Convex Geometry discovered by V. D. Milman in the mid 1980s 
(cf. \cite{bm-M1, bm-M2, bm-M3, bm-Pis}). It has connections with high dimensional phenomena, 
and we refer an interested reader to \cite{bm-BKM, bm-KT, bm-KM, bm-AMO}. 
The questions concerning possible description of the maps $T_1$ and $T_2$ 
and related isotropic properties of the normalized Gaussian measures are 
discussed in \cite{bm-Bob2}.

Based on \eqref{bm-rBM}, and involving Berwald's inequality in the form of 
C. Borell \cite{bm-Bor1}, the inverse entropy power inequality \eqref{bm-repi-q} has been 
established recently \cite{bm-BM1, bm-BM3} for the class of all probability 
distributions having log-concave densities. Involving additionally a general 
submodularity property of entropy \cite{bm-Mad}, it turned out also possible 
to consider more general densities of the form
\begin{equation}\label{bm-cvx}
f(x) = V(x)^{-\beta}, \qquad x \in \RR^n,
\end{equation}
where $V$ are positive convex functions on $\RR^n$ and $\beta \geq n$ is a given
parameter. More precisely, the following statement can be found in \cite{bm-BM3}.

\begin{theorem}\label{bm-thm:repi}
Let $X$ and $Y$ be independent random vectors 
in $\RR^n$ with densities of the form $\eqref{bm-cvx}$ with $\beta \geq 2n+1$, 
$\beta \geq \beta_0 n$ $(\beta_0 > 2)$. There exist linear volume preserving 
maps $T_i:\RR^n \rightarrow \RR^n$ such that 
\begin{equation}\label{bm-repi}
H\big(\widetilde X + \widetilde Y\big)\, \leq\, C_{\beta_0}\, (H(X) + H(Y)),
\end{equation}
where $\widetilde X = T_1(X)$, $\widetilde Y = T_2(Y)$, and where
$C_{\beta_0}$ is a constant, depending on $\beta_0$, only.
\end{theorem}

The question of what maps $T_1$ and $T_2$ can be used in Theorem~\ref{bm-thm:repi} is 
rather interesting, but certainly the maps that put the distributions 
of $X$ and $Y$ in $M$-position suffice (see \cite{bm-BM3} for terminology and 
discussion). In a more relaxed form, one needs to have in some sense 
``similar'' positions for both distributions. For example, when considering 
identically distributed random vectors, there is no need to appeal 
in Theorem~\ref{bm-thm:repi} to some (not very well understood) affine volume-preserving
transformations, since the distributions of $X$ and $Y$ have the same 
$M$-ellipsoid. In other words, we have for $X$ and $Y$ drawn independently 
from the {\it same} distribution (under the same assumption on form of 
density as Theorem~\ref{bm-thm:repi}) that
\begin{equation}\label{bm-repi-sum}
H(X + Y)\, \leq\, C_{\beta_0}\, (H(X) + H(Y))= 2 C_{\beta_0}\, H(X).
\end{equation}
Since the distributions of $X$ and $-Y$ also have the same $M$-ellipsoid,
it is also true that
\begin{equation}\label{bm-repi-diff}
H(X - Y)\, \leq\, C_{\beta_0}\, (H(X) + H(Y)) = 2 C_{\beta_0}\, H(X).
\end{equation}

We strengthen this observation by providing a quantitative version 
with explicit constants below (under, however, a convexity condition on 
the convolved measure). Moreover, one can give a short and relatively 
elementary proof of it without appealing to Theorem~\ref{bm-thm:repi}.

\begin{theorem}\label{bm-thm:sumdiff}
Let $X$ and $Y$ be independent identically
distributed random vectors in $\RR^n$ with finite entropy. Suppose that 
$X-Y$ has a probability density function of the form $\eqref{bm-cvx}$ with
$\beta\geq \max\{n+1,\beta_0 n\}$ for some fixed $\beta_0>1$. Then
$$
H(X-Y)\leq D_{\beta_0} H(X) 
$$
and
$$
H(X+Y)\leq D_{\beta_0}^2 H(X) ,
$$
where
$D_{\beta_0}=\exp(\frac{2\beta_0}{\beta_0-1})$.
\end{theorem}

Let us return to Theorem~\ref{bm-thm:repi} and the class of distributions involved 
there. For growing $\beta$, the families \eqref{bm-cvx} shrink and converge in the 
limit as $\beta \rightarrow +\infty$ to the family of log-concave densities 
which correspond to the class of log-concave probability measures.
Through inequalities of the Brunn-Minkowski-type, the latter class was 
introduced by A. Pr\'ekopa, while the general case $\beta \geq n$
was studied by C. Borell \cite{bm-Bor2, bm-Bor3}, cf. also \cite{bm-BL, bm-Bob1}.
In \cite{bm-Bor2, bm-Bor3} it was shown that probability
measures $\mu$ on $\RR^n$ with densities \eqref{bm-cvx} (and only they, once $\mu$ 
is absolutely continuous) satisfy the geometric inequality
\begin{equation}\label{bm-cvx-meas}
\mu \big (tA + (1-t)B \big ) \geq
\big [\,t\mu(A)^\kappa + (1-t)\mu(B)^\kappa\big ]^{1/\kappa}
\end{equation}
for all $t \in (0,1)$ and for all Borel measurable sets $A,B \subset \RR^n$,
with negative power
$$
\kappa = - \frac{1}{\beta - n}.
$$
Such $\mu$'s form the class of so-called $\kappa$-concave measures.
In this hierarchy the limit case $\beta = n$ corresponds to $\kappa = -\infty$
and describes the largest class of measures on $\RR^n$, called {\it convex}, 
in which case \eqref{bm-cvx-meas} turns into
$$
\mu(tA + (1-t)B) \geq \min \{\mu(A),\mu(B)\}.
$$
This inequality is often viewed as the weakest convexity hypothesis
about a given measure $\mu$.

One may naturally wonder whether or not it is possible to relax the
assumption on the range of $\beta$ in \eqref{bm-repi}-\eqref{bm-repi-diff}, or even to remove 
any convexity  hypotheses. In this note we show that this is impossible
already for the class of all one-dimensional convex probability distributions.
Note that in dimension one there are only two admissible linear transformations, 
$\widetilde X = X$ and $\widetilde X = -X$, so that one just wants 
to estimate $H(X+Y)$ or $H(X-Y)$ from above in terms of $H(X)$. 
As a result, the following statement demonstrates that 
Theorem~\ref{bm-thm:repi} and its particular cases \eqref{bm-repi-sum}-\eqref{bm-repi-diff} are false over the full 
class of convex measures.

\begin{theorem}\label{bm-thm:neg}
For any constant $C$, there is a convex probability 
distribution $\mu$ on the real line with a finite entropy, such that
$$
\min\{H(X+Y), H(X-Y)\} \geq C\,H(X),
$$
where $X$ and $Y$ are independent random variables, distributed according
to $\mu$.
\end{theorem}

A main reason for $H(X+Y)$ and $H(X-Y)$ to be much larger than $H(X)$
is that the distributions of the sum $X+Y$ and the difference $X-Y$ may
lose convexity properties, when the distribution $\mu$ of $X$
is not ``sufficiently convex''. For example, in terms of the convexity
parameter $\kappa$ (instead of $\beta$), the hypothesis of Theorem~\ref{bm-thm:repi}
is equivalent to
$$
\kappa \geq - \frac{1}{(\beta_0 - 1)n} \ \ (\beta_0 > 2), \qquad
\kappa \geq - \frac{1}{n+1}.
$$
That is, for growing dimension $n$ we require that $\kappa$ be sufficiently close 
to zero (or the distributions of $X$ and $Y$ should be close to the class
of log-concave measures). These conditions ensure that
the convolution of $\mu$ with the uniform distribution on a proper
(specific) ellipsoid remains to be convex, and its convexity parameter
can be controled in terms of $\beta_0$ 
(a fact used in the proof of Theorem~\ref{bm-thm:repi}). 
However, even if $\kappa$ is close to zero, one cannot guarantee that $X + Y$ or 
$X - Y$ would have convex distributions.

We prove Theorem~\ref{bm-thm:sumdiff} in Section~\ref{bm-sec:RS} and Theorem~\ref{bm-thm:neg} in Section~\ref{bm-sec:pf}, and then 
conclude in Section~\ref{bm-sec:rmk} with remarks on the relationship between Theorem~\ref{bm-thm:neg} 
and recent results about Cramer's characterization of the normal law.

\section{A ``difference measure'' inequality for convex measures}
\label{bm-sec:RS}

Given two convex bodies $A$ and $B$ in $\RR^n$, introduce 
$A-B = \{x-y: x \in A, \ y \in B\}$. 
In particular, $A-A$ is called the "difference body" of $A$.
Note it is always symmetric about the origin.

The Rogers-Shephard inequality \cite{bm-RS} states that, for any convex body 
$A\subset \RR^n$,
\begin{equation}\label{bm-rs}
|A - A|\, \leq\, C_{2n}^n\, |A|,
\end{equation}
where $C_n^k = \frac{n!}{k! (n-k)!}$ denote usual combinatorial coefficients. 
Observe that putting the Brunn-Minkowski inequality and \eqref{bm-rs} together 
immediately yields that
\begin{eqnarray*}
2 \leq \frac{|A-A|^{\frac{1}{n}}}{|A|^{\frac{1}{n}}} \leq \left[C_{2n}^n\right]^\frac{1}{n} < 4,
\end{eqnarray*}
which constrains severely the volume radius of the difference body of $A$
relative to that of $A$ itself. In analogy to the Rogers-Shephard inequality,
we ask the following question for entropy of convex measures.

\vskip3mm
\noindent{\bf Question.} 
{\it Let $X$ and $Y$ be independent random vectors in $\RR^n$, which are 
identically distributed with density
$V^{-\beta}$, with $V$ positive convex, and $\beta \geq n+\gamma$. 
For what range of $\gamma>0$ is it true that
$H(X-Y) \leq C_{\gamma} H(X)$, for some constant $C_{\gamma}$ 
depending only on $\gamma$}?
\vskip3mm

Theorems~\ref{bm-thm:sumdiff} and \ref{bm-thm:neg} partially answer this question. To prove the former, 
we need the following lemma about convex measures, proved in \cite{bm-BM2}.

\begin{lemma}\label{bm-lem:cvx-maxnorm}
Fix $\beta_0 > 1$. Assume a random vector $X$ in $\RR^n$ has 
a density $f = V^{-\beta}$, where $V$ is a positive convex function
on the supporting set.
If $\beta \geq n+1$ and $\beta \geq \beta_0 n$, then
\begin{equation}\label{bm-maxnorm}
\log\, \|f\|_{\infty}^{-1} \leq\,
h(X)  \leq\, c_{\beta_0} n + \log\, \|f\|_{\infty}^{-1} ,
\end{equation}
where one can take for the constant $c_{\beta_0}=\frac{\beta_0}{\beta_0-1}$.
\end{lemma}

In other words, for sufficiently convex probability measures, the entropy may 
be related to the $L^\infty$-norm $\|f\|_{\infty} = \sup_x f(x)$ of the density 
$f$ (which is necessarily finite). Observe that the left inequality in \eqref{bm-maxnorm} 
is general: It trivially holds without any convexity assumption. On the other 
hand, the right inequality is  an asymptotic version of a result from \cite{bm-BM2}
about extremal role of the multidimensional Pareto distributions.

Now, let $f$ denote the density of the random variable
$W = X-Y$ in Theorem~\ref{bm-thm:sumdiff}. It is symmetric (even) and thus
maximized at zero, by the convexity hypothesis. Hence, by Lemma~\ref{bm-lem:cvx-maxnorm},
\begin{eqnarray*}
h(W)\leq \log \|f\|_{\infty}^{-1} + c_{\beta_0} n
= \log f(0)^{-1} + c_{\beta_0} n.
\end{eqnarray*}
But, if $p$ is the density of $X$, then $f(0)= \int_{\RR^n} p(x)^2\, dx$,
and hence
\begin{eqnarray*}
\log f(0)^{-1} = -\log \int_{\RR^n} p(x) \cdot p(x)\, dx
\leq \int_{\RR^n} p(x) [-\log p(x)]\, dx 
\end{eqnarray*}
by using Jensen's inequality. Combining the above two
displays immediately yields the first part of Theorem~\ref{bm-thm:sumdiff}.

To obtain the second part, we need an observation from \cite{bm-MK} that follows from 
the following lemma on the submodularity of the entropy of sums proved in \cite{bm-Mad}.

\begin{lemma}\label{bm-lem:submod}
Given independent random vectors $X,Y,Z$
in $\RR^n$ with absolutely continuous distributions, we have
$$
h(X+Y+Z) + h(Z) \leq h(X+Z) + h(Y+Z),
$$
provided that all entropies are well-defined and finite.
\end{lemma}

Taking $X$, $Y$ and $-Z$ to be identically distributed, and using the
monotonicity of entropy (after adding an independent summand), we obtain
\begin{eqnarray*}
h(X+Y)+h(Z)\leq h(X+Y+Z) + h(Z) \leq h(X+Z) + h(Y+Z)
\end{eqnarray*}
and hence
\begin{eqnarray*}
h(X+Y)+h(X) \leq 2 h(X-Y).
\end{eqnarray*}
(This is the relevant observation from \cite{bm-MK}.) Combining this bound with 
the first part of Theorem~\ref{bm-thm:sumdiff} immediately gives the second part.


It would be more natural to state Theorem~\ref{bm-thm:sumdiff} under a shape condition on 
the distribution of $X$ rather than on that of $X-Y$, but for this we need 
to have better understanding of the convexity parameter of the convolution 
of two $\kappa$-concave measures when $\kappa<0$. 

Observe that in the log-concave case of Theorem~\ref{bm-thm:sumdiff} (which is the case of 
$\beta\rightarrow\infty$, but can easily be directly derived in the same way without 
taking a limit), one can impose only a condition on the distribution of $X$ 
(rather than that of $X-Y$) since closedness under convolution is guaranteed 
by the Pr\'ekopa-Leindler inequality.

\begin{corollary}\label{bm-cor:ent-rs}
Let $X$ and $Y$ be independent random vectors in 
$\RR^n$ with log-concave densities. Then
\begin{eqnarray*}
h(X-Y) & \leq & h(X) + n, \\
h(X+Y) & \leq & h(X) + 2n.
\end{eqnarray*}
\end{corollary}

In particular, observe that putting the entropy power inequality \eqref{bm-epi} 
and Corollary~\ref{bm-cor:ent-rs} together immediately yields that
\begin{eqnarray*}
2 \leq \frac{H(X-Y)}{H(X)} \leq e^2, 
\end{eqnarray*}
which constrains severely the entropy power of the ``difference measure'' 
of $\mu$ relative to that of $\mu$ itself.

\section{Proof of Theorem~\ref{bm-thm:neg}}
\label{bm-sec:pf}

Given a (large) parameter $b>1$, let a random variable $X_b$ have
a truncated Pareto distribution $\mu$, namely, with the density
$$
f(x) = \frac{1}{x\log b}\ 1_{\{1 < x < b\}}(x).
$$
By the construction, $\mu$ is supported on a bounded interval $(1,b)$ and 
is convex.

First we are going to test the inequality 
\begin{equation}\label{bm-testsum}
H(X_b + Y_b) \leq CH(X_b)
\end{equation}
for growing $b$, where $Y_b$ is an independent copy of $X_b$. Note that
\begin{eqnarray*}
h(X_b) & = & \int_1^b f(x)\, \log(x \log b)\,dx \\
  & = & 
\log \log b + \frac{1}{\log b} \int_1^b \frac{\log x}{x}\,dx
  \ = \
\log \log b + \frac{1}{2}\ \log b,
\end{eqnarray*}
so $H(X_b) = b\, \log^2 b$.

Now, let us compute the convolution of $f$ with itself. The sum $X_b + Y_b$ 
takes values in the interval $(2,2b)$. Given $2 < x < 2b$, we have
$$
g(x) = (f*f)(x) = \int_{-\infty}^{+\infty} f(x-y)f(y)\,dy =
\frac{1}{\log^2 b} \int_\alpha^\beta \frac{dy}{(x-y)y},
$$
where the limits of integration are determined to satisfy the constraints
$1 < y < b$, $1 < x-y < b$. So, 
$$
\alpha = \max(1,x-b), \quad \beta = \min(b,x-1),
$$
and using $\frac{1}{(x-y)y} = \frac{1}{x}\,(\frac{1}{y} + \frac{1}{x-y})$, we
find that
\begin{eqnarray*}
g(x) & = & 
\frac{1}{x\log^2 b}\ \big(\log(y) - \log(x-y)\big|_{x = \alpha}^\beta
 \ = \ \frac{1}{x\log^2 b}\ \log\frac{y}{x-y}\,\bigg|_{x = \alpha}^\beta \\
 & = & 
\frac{1}{x\log^2 b}\ 
\bigg(\log\frac{\beta}{x-\beta} - \log\frac{\alpha}{x-\alpha}\bigg). \\ 
\end{eqnarray*}
Note that 
$x - \alpha = x - \max(1,x-b) = \min(b,x-1) = \beta$.
Hence,
$$
g(x) = \frac{2}{x\log^2 b}\ \log\frac{\beta}{\alpha} =
\frac{2}{x\log^2 b}\ \log\frac{\min(b,x-1)}{\max(1,x-b)}.
$$
Equivalently,
\begin{eqnarray*}
g(x) = \frac{2}{x\log^2 b}\ \log(x-1), & {\rm for} & 2 < x < b+1, \\
g(x)\, =\, \frac{2}{x\log^2 b}\ \log \frac{b}{x-b}, & {\rm for} & b+1 < x < 2b.
\end{eqnarray*}

Now, on the second interval $b+1 < x < 2b$, we have
$$
g(x) \leq \frac{2}{x\log^2 b}\, \log b = \frac{2}{x\log b} < \frac{2}{(b+1)\log b} < 1,
$$
where the last bound holds for $b \geq e$, for example. Similarly,
on the first interval $2 < x < b+1$, using $\log(x-1) < \log b$, we get
$$
g(x) \leq \frac{2}{x\log b}  < \frac{1}{\log b} \leq 1.
$$
Thus, as soon as $b \geq e$, we have $g \leq 1$ on the support interval. From this,
$$
h(X_b + Y_b) = \int_2^{2b} g(x) \log(1/g(x))\,dx \geq
\int_2^{b} g(x) \log(1/g(x))\,dx.
$$
Next, using on the first interval the bound
$g(x) \leq \frac{2}{x\log b} \leq \frac{1}{x}$, valid for $b \geq e^2$, we get
for such values of $b$ that
$$
h(X_b + Y_b) \geq \int_2^{b} g(x) \log x\,dx = \frac{2}{\log^2 b}
\int_2^{b} \frac{\log(x-1)\, \log x}{x}\,dx.
$$
To further simplify, we may write $x - 1 \geq \frac{x}{2}$, which gives
\begin{eqnarray*}
\int_2^{b} \frac{\log(x-1)\, \log x}{x}\,dx
 & \geq &
\int_2^{b} \frac{\log^2 x}{x}\,dx - \log 2 \int_2^{b} \frac{\log x}{x}\,dx \\
 & = &
\frac{1}{3}\, \big(\log^3 b - \log^3 2\big) -
\frac{\log 2}{2}\ \big(\log^2 b - \log^2 2\big) \\
 & > &
\frac{1}{3}\, \log^3 b - \frac{\log 2}{2}\, \log^2 b.
\end{eqnarray*}
Hence,
$
h(X_b + Y_b) > \frac{2}{3}\, \log b - \log 2,
$
and so
$$
H(X_b + Y_b) > \frac{1}{2}\,b^{4/3} \quad (b \geq e^2).
$$
In particular,
$$
\frac{H(X_b + Y_b)}{H(X_b)} > \frac{b^{1/3}}{2\log b} \rightarrow +\infty, \quad
{\rm as} \ \ b \rightarrow +\infty.
$$
Hence, the inequality \eqref{bm-testsum} may not hold for large $b$ with any 
prescribed value of $C$.

To test the second bound
\begin{equation}\label{bm-testdiff}
H(X_b - Y_b) \leq CH(X_b),
\end{equation}
one may use the previous construction. The random variable 
$X_b - Y_b$ can take any value in the interval $|x| < b-1$, where it is described by 
the density
$$
h(x) = \int_{-\infty}^{+\infty} f(x+y)f(y)\,dy =
\frac{1}{\log^2 b} \int_\alpha^\beta \frac{dy}{(x+y)y}.
$$
Here the limits of integration are determined to satisfy
$1 < y < b$ and $1 < x+y < b$. So, assuming for simplicity that $0<x< b-1$, the
limits are
$$
\alpha = 1, \quad \beta = b-x.
$$
Writing $\frac{1}{(x+y)y} = \frac{1}{x}\,(\frac{1}{y} - \frac{1}{x+y})$, we
find that
$$
h(x) =  \frac{1}{x\log^2 b}\ \big(\log(y) - \log(x+y)\big|_{x = \alpha}^\beta = 
\frac{1}{x\log^2 b}\ \log\frac{(b-x) (x+1)}{b}. \\ 
$$
It should also be clear that
$$
h(0) =
\frac{1}{\log^2 b} \int_1^{b} \frac{dy}{y^2} = \frac{1 - \frac{1}{b}}{\log^2 b}.
$$

Using $\log \frac{(b-x) (x+1)}{b} < \log(x+1) < x$, we obtain that
$h(x) < \frac{1}{\log^2 b} \leq 1$, for $b \geq e^2$.

In this range, since $\frac{(b-x) (x+1)}{b} < b$, we also have that 
$h(x) \leq \frac{1}{x\log b} \leq \frac{1}{x}$. Hence, in view of the
symmetry of the distribution of $X_b - Y_b$,
\begin{eqnarray*}
h(X_b - Y_b) & = & 2 \int_0^{b-1} h(x) \log(1/h(x))\,dx \\
 & \geq &
2\int_0^{b/2} h(x) \log x\,dx \\ 
 & = & 
 \frac{2}{\log^2 b} \int_2^{b/2} \frac{\log x}{x}\,\ \log\frac{(b-x) (x+1)}{b}\,dx.
\end{eqnarray*}
But for $0 < x < b/2$, 
$$
\log\frac{(b-x) (x+1)}{b} > \log\frac{x+1}{2} > \log x - \log 2,
$$
so
\begin{eqnarray*}
h(X_b - Y_b) & > & 
\frac{2}{\log^2 b} \int_2^{b/2} \frac{\log^2 x - \log 2\, \log x}{x}\,dx \\
 & = &
\frac{2}{\log^2 b}\, \bigg(\frac{1}{3}\, (\log^3 (b/2) - \log^3 2) -
\frac{\log 2}{2}\, (\log^2 (b/2) - \log^2 2)\bigg) \\
 & > &
\frac{2}{\log^2 b}\, \bigg(\frac{1}{3}\,\log^3 (b/2) -
\frac{1}{2}\,\log^2 (b/2)\bigg) \\
 & \sim &
\frac{2}{3}\, \log b.
\end{eqnarray*}
Therefore, like on the previous step, $H(X_b - Y_b)$ is bounded from below 
by a function, which is equivalent to $b^{4/3}$. Thus, for large $b$, 
the inequality \eqref{bm-testdiff} may not hold either.

Theorem~\ref{bm-thm:neg} is proved.

\section{Remarks}
\label{bm-sec:rmk}

For a random variable $X$ having a density, consider 
the entropic distance from the distribution of $X$ to normality
$$
D(X) = h(Z) - h(X),
$$
where $Z$ is a normal random variable with parameters $\EE Z = \EE X$, 
${\rm Var}(Z) = {\rm Var}(X)$. This functional is well-defined for the class
of all probability distributions on the line with finite second moment, 
and in general $0 \leq D(X) \leq +\infty$.

The entropy power inequality implies that
\begin{eqnarray}\label{bm-epi-D}
D(X+Y) & \leq & 
\frac{\sigma_1^2}{\sigma_1^2 + \sigma_2^2}\, D(X) +
\frac{\sigma_2^2}{\sigma_1^2 + \sigma_2^2}\, D(X) \nonumber \\
 & \leq &
\max(D(X),D(Y)),
\end{eqnarray}
where $\sigma_1^2 = {\rm Var}(X)$, $\sigma_2^2 = {\rm Var}(Y)$.

In turn, if $X$ and $Y$ are identically distributed, then Theorem~\ref{bm-thm:neg}
reads as follows: For any positive constant $c$, there exists a convex 
probability measure $\mu$ on $\RR$ with $X,Y$ independently distributed 
according to $\mu$, with
$$
D(X \pm Y) \leq D(X) - c.
$$
This may be viewed as a strengthened variant of \eqref{bm-epi-D}.
That is, in Theorem~\ref{bm-thm:neg} we needed to show that both $D(X + Y)$ and $D(X - Y)$
may be much smaller than $D(X)$ in the additive sense. In particular, 
$D(X)$ has to be very large when $c$ is large.
For example, in our construction of the previous section
$$
\EE X_b = \frac{b-1}{\log b}, \quad \EE X_b^2 = \frac{b^2-1}{2 \log b},
$$
which yields
$$
D(X_b) \sim \frac{3}{2}\, \log b, \quad D(X_b + Y_b) \sim \frac{4}{3}\, \log b,
$$
as $b \rightarrow +\infty$.

In \cite{bm-BCG1, bm-BCG2} a slightly different question, raised by M. Kac and H.~P.~McKean 
\cite{bm-McK} (with the desire to quantify in terms of entropy the Cramer 
characterization of the normal law), has been answered. Namely, it was shown 
that $D(X+Y)$ may be as small as we wish, while $D(X)$ is 
separated from zero. In the examples of \cite{bm-BCG2}, $D(X)$ is
of order 1, while for Theorem~\ref{bm-thm:neg} it was necessary to use large values for 
$D(X)$, arbitrarily close to infinity.
In addition, the distributions in \cite{bm-BCG1, bm-BCG2} are not convex.


\end{document}